\documentclass[hidelinks]{siamart190516}
\usepackage{graphicx,amsfonts,upquote,verbatim}

\def\qed{~\hbox{\vrule width 3pt depth 4 pt height 2.5 pt}}
\def\pput(#1,#2)#3{\noindent\smash{\raise#2pt\hbox to 0pt
   {\kern #1pt #3\hss}}\ignorespaces}

\usepackage{tikz} 
\usepackage[ruled,vlined]{algorithm2e}
\newcounter{example}[section]

\newcommand{\reals}{\mathbb{R}}
\newcommand{\complexes}{\mathbb{C}}

\def\Re{\hbox{\rm Re\kern .3pt}}
\def\Im{\hbox{\rm Im}}
\def\gradient{\boldsymbol{\nabla}\kern -.7pt}
\def\laplacian{\Delta\kern -.3pt}

\title{Lightning Stokes solver}

\author{Pablo D. Brubeck\thanks{\texttt{brubeckmarti@maths.ox.ac.uk},
Mathematical Institute, University of Oxford, Oxford OX2 6GG, UK.}
\and
Lloyd N.~Trefethen\thanks{\texttt{trefethen@maths.ox.ac.uk},
Mathematical Institute, University of Oxford, Oxford OX2 6GG, UK.}}

\begin{document}

\maketitle

\begin{abstract}
Gopal and Trefethen recently introduced ``lightning solvers''
for the 2D Laplace and Helmholtz equations, based on rational
functions with poles exponentially clustered near singular corners.
Making use of the Goursat representation in terms of analytic
functions, we extend these methods to the biharmonic equation,
specifically to 2D Stokes flow.  Solutions to model problems
are computed to 10-digit accuracy in less than a second of laptop time.
As an illustration of the high accuracy, we resolve two or more
counter-rotating Moffatt eddies near a singular corner.
\end{abstract}


\pagestyle{myheadings}
\thispagestyle{plain}
\markboth{\sc Brubeck and Trefethen}
{\sc Lightning Stokes solver}

\section{\label{secintro}Introduction}

Consider the biharmonic equation in a
two-dimensional domain $\Omega\subseteq\reals^2$,
\begin{equation} \label{eq:bih}
\laplacian^2 \psi = \psi_{xxxx} +2\kern .3pt\psi_{xxyy} +\psi_{yyyy} = 0,
\end{equation}
where the subscripts denote partial differentiation.
Since the equation is of fourth order, two boundary conditions are
imposed at each point of $\partial\kern.5pt\Omega$. 
On one (nonempty) part $\Gamma_1\subset \partial\kern.5pt\Omega$,
the value of the function and one component of the gradient are specified,
and on the remainder $\Gamma_2 = \partial\kern.5pt\Omega \backslash \Gamma_1$,
both components of the gradient are specified:
\begin{equation}\label{eq:bcsa}
\psi = h(x,y), \quad \textbf a(x,y) \cdot \gradient\psi = k(x,y),\quad (x,y) \in\Gamma_1,\\
\end{equation}
\begin{equation}\label{eq:bcsb}
\gradient \psi = \textbf g (x,y), \quad (x,y)\in\Gamma_2
\end{equation}
(boldface denotes vectors).  Functions that satisfy (\ref{eq:bih})
are called {\em biharmonic\/} and are smooth in the interior
of~$\Omega$, but point singularities may arise on the boundary
$\partial\kern.5pt\Omega$ when it contains corners, or when the
data $h$, $k$, $\mathbf{a}$, $\textbf g$ have singularities,
or at junctions between $\Gamma_1$ and $\Gamma_2$.

In this paper we propose a numerical method to solve
(\ref{eq:bih})--(\ref{eq:bcsb}) by generalizing the recently
introduced ``lightning solvers'' for the 2D Laplace and Helmholtz
equations \cite{gopal19new,gopal19,clustering}. The main advantage
of this class of methods is that they can handle domains with
corners without requiring any detailed analysis, and still achieve
high accuracy, with root-exponential convergence as a function
of the number of degrees of freedom.

The work is structured as follows. In section~\ref{sec:goursat},
we review the Goursat representation, which allows one to write
a biharmonic function in terms of two analytic functions when
$\Omega$ is simply connected.  In section~\ref{sec:physics} we
present the application to 2D Stokes flow, where it is the stream
function that is biharmonic.   The biharmonic equation also arises
in linear elasticity problems, but elasticity problems are not
considered in this paper~\cite{goodier}.

Our method consists of approximating both Goursat functions
by rational functions with fixed poles and finding
coefficients that best satisfy the boundary conditions
(\ref{eq:bcsa})--(\ref{eq:bcsb}) in a least-squares
sense.  A detailed description of the method is given in
section~\ref{sec:method}, and numerical results for several
cavity flows are presented in section~\ref{sec:results}. The
method gets high accuracy for simple geometries with great speed,
and we illustrate this by showing its ability to resolve two or
more Moffatt eddies in corners.  Section~\ref{sec:analytic}
examines the behavior of computed solutions outside
the domain $\Omega$, where they are also defined, and 
Section~\ref{sec:unbounded} presents further examples involving
channel flows.
The strengths and weaknesses of the method are considered in
the discussion section, section~\ref{sec:conclusion}, along with its prospects for extension
to multiply-connected and other geometries.
The appendix lists a template program for these calculations which 
readers can adapt to their own problems.

The Goursat representation has been applied for numerical solution of
biharmonic problems by Luca and coauthors~\cite{luca18,lucals} and Kazakova
and Petrov~\cite{kaz}, among others.  However, it is used far
less than more general tools such as the finite element method
and integral equations.

\section{Biharmonic functions in simply-connected domains\label{sec:goursat}}

The key to extending the lightning solver to biharmonic
problems is the Goursat representation, introduced by the French
mathematician \'Edouard Goursat, author of a celebrated {\em Cours
d'analyse math\'ematique}~\cite{goursat98}.  This allows one to
represent a biharmonic function in terms of two analytic functions.
Throughout the rest of the paper, we combine $x$ and $y$ in the
complex variable $z = x + iy$.  We begin with a precise statement
of the Goursat representation in the form of a theorem, along with
two proofs, which are helpful in providing different perspectives
on what is going on.  The first proof was published by Goursat
in a note of barely more than one page~\cite{goursat98}, and the
second was given by Muskhelishvili in~\cite{musk19}.  Both can
also be found in textbooks such as \cite{carrier05,musk77}.
See also chapter~7 of~\cite{langdev}, where, however, the name
Goursat is not mentioned.  Although we label them as ``proofs,''
we do not spell out the arguments in detail.  For a more rigorous
discussion, see~\cite{musk77}.  Our use of the imaginary as
opposed to the real part in (\ref{eq:goursat}) is arbitrary.

\begin{theorem}[Goursat representation]
\label{th:goursat}
Let $\Omega\subseteq\complexes$ be a simply-connected open region, and let
$\psi : \Omega \to \reals$ be a biharmonic function. Then $\psi$ can
be represented in terms of complex
functions $f(z)$ and $g(z)$ analytic in $\Omega$ by the formula
\begin{equation}\label{eq:goursat}
\psi(x,y) = \psi(z,\overline{z}\kern .5pt ) =
\Im \kern -.8pt \left\{ \overline{z}f(z) + g(z)\right\}.
\end{equation}
The functions $f(z)$ and $g(z)$ are uniquely defined up to addition of
arbitrary terms $\gamma z+C$ to $f(z)$ and $\overline{C}z+\alpha$ to $g(z)$,
where $\gamma,\alpha\in\reals$ and\/ $C\in\complexes$.
\end{theorem}

\medskip

{\em Goursat's proof}~\cite{goursat98}.
First, we note that the biharmonic operator can be written in terms of
Wirtinger derivatives,
\begin{displaymath}
\laplacian^2 \psi = 16\kern .7pt \frac{\partial^4 \psi}{\partial\overline{z}^2 \partial z^2}.
\end{displaymath}
Treating $z$ and $\overline{z}$ as independent, one can solve this
equation by separation of variables. 
The ansatz $\psi = A(z)B(\overline{z})$ leads to the equation
\begin{displaymath}
\frac{\partial^4 \psi}{\partial\overline{z}^2 \partial z^2} 
   = A''(z) B''(\overline{z}) = 0,
\end{displaymath} 
which implies that either $A''(z)=0$ or $B''(\overline{z})=0$. The solution
spaces for these two cases are spanned by $\{1,z\}$ for $A''(z)=0$  and
$\{1,\overline{z}\kern .3pt \}$ for $B''(\overline{z})=0$,
so the general solution takes the form
\begin{displaymath}
\psi(z,\overline{z}\kern .3pt )
=  A_1(z) + \overline{z} A_2(z) + B_3(\overline{z}) + z B_4(\overline{z}),
\end{displaymath}
where $A_1(z), A_2(z), B_3(\overline{z}), B_4(\overline{z})$ are arbitrary
analytic functions.
Since $\psi$ must be real, i.e., $\psi=\overline{\psi}$, we must have
\begin{displaymath}
\psi(z,\overline{z}\kern .3pt)
= A_1(z)  + \overline{z} A_2(z)  + \overline{A_1(z)} + z \overline{A_2(z)} 
                 = 2\kern .5pt \Re\kern -.8pt \left\{A_1(z) + \overline{z} A_2(z)\right\},
\end{displaymath}
which implies (\ref{eq:goursat}) by letting $f(z) = -2\kern .3pt i A_2(z)$
and $g(z) = -2\kern .3pt i A_1(z)$. We finally note that, when $\psi$ is given, $f(z)$
and $g(z)$ are determined up to additive linear terms. To be precise,
the substitutions $f(z)\to f(z)+\gamma z + C$ together with $g(z)\to
g(z)+\overline{C}z+\alpha$ for $\gamma,\alpha\in\reals$ and $C\in\complexes$
leave $\psi$ invariant.
\qed

\medskip

{\em Muskhelishvili's proof} \cite{musk19}.
If $\omega=-\laplacian \psi$, then
$-\laplacian\kern .8pt \omega = \laplacian^2\psi = 0$. 
Define $p$ as the harmonic conjugate of $\omega$,
determined up to an arbitrary constant $p_0\in\reals$ and
satisfying the Cauchy--Riemann conditions
\begin{displaymath}
\frac{\partial\kern.3pt\omega}{\partial x} = \frac{\partial p}{\partial y}, \quad
\frac{\partial\kern.3pt\omega}{\partial y} = -\frac{\partial p}{\partial x}\kern 1pt .
\end{displaymath}
Then the expression
$h(z) = p(x,y)-i\kern .4pt\omega(x,y)$
represents an analytic function of $z=x+iy$.  Put
\begin{displaymath}
f(z) = \frac{1}{4} \int_a^z h(w) dw = f_1 + if_2,
\end{displaymath}
where $a\in\Omega$ is arbitrary, causing $f(z)$ to be defined up to a term
$\gamma z +C$ where $4\gamma=p_0$ and $C\in\complexes$. Then
\begin{displaymath}
f'(z) = \frac{\partial f_1}{\partial x} + i \kern .3pt\frac{\partial f_2}{\partial x} 
   = -i\kern .3pt \frac{\partial f_1}{\partial y} + \frac{\partial f_2}{\partial y} 
   = \frac{1}{4} \left(\kern .7pt p-i\kern .4pt\omega\right),
\end{displaymath}
from which we calculate, using the fact that $\laplacian f_1=\laplacian f_2=0$,
\begin{displaymath}
\laplacian\left(y f_1\right) = 2\kern .5pt\frac{\partial f_1}{\partial y} 
   = \frac{\omega}{2},\quad \laplacian\kern -1pt \left(x f_2\right) 
   = 2\kern .5pt\frac{\partial f_2}{\partial x} = -\frac{\omega}{2}.
\end{displaymath}
Thus the function $\psi + yf_1 - xf_2$ is harmonic, since
\begin{displaymath}
\laplacian \left(\psi + yf_1 - xf_2\right) 
   = -\omega + \frac{\omega}{2} + \frac{\omega}{2} =0.
\end{displaymath}
If we write $\psi + yf_1 - xf_2$ as $\Im \kern .8pt g(z)$, where
$g(z)$ is analytic, we obtain as required
\begin{displaymath}
\psi = x f_2(z) -y f_1(z) + \Im \kern .8pt g(z)
     = \Im\kern -.9pt \left\{\overline{z}f(z) + g(z)\right\}. \qed
\end{displaymath}

\section{Stokes flow and Moffatt eddies\label{sec:physics}}

The application of the biharmonic
equation and the Goursat representation to Stokes flow is
an old idea and is mentioned for example as an exercise in~\cite{ockendon95}.
Let $\textbf u=(u,v)^T$ be the velocity field of a steady
incompressible 2D fluid flow, and let $p$ be the associated pressure field,
which for mathematical purposes
we may take to be defined up to a constant.
In the limit of zero Reynolds number and in
the absence of body forces, the steady-state equations for $\textbf u$ and $p$ are
\begin{equation}
-\laplacian \textbf u + \gradient p = \textbf 0,\quad
\gradient\cdot\textbf u = 0, \label{eq:incomp}
\end{equation}
representing conservation of momentum and incompressibility, respectively.
We now introduce the stream function $\psi$ by defining 
\begin{equation}\label{eq:vel}
u=\frac{\partial\psi}{\partial y}, \quad 
v=-\frac{\partial\psi}{\partial x},
\end{equation}
or equivalently,
\begin{equation}\label{eq:vel2}
\gradient \psi = \begin{pmatrix} -v \\ u \end{pmatrix}.
\end{equation}
Note that $\psi$ is unique up to an additive constant, which may be
chosen arbitrarily.
It is easy to verify that the second equation
of (\ref{eq:incomp}) is automatically satisfied.
Now define the vorticity as
\begin{equation} \label{eq:poisson}
\omega = \frac{\partial v}{\partial x}-\frac{\partial u}{\partial y}  = -\laplacian \psi.
\end{equation}
From the first equation of (\ref{eq:incomp}) we have $-\laplacian u = -\partial p/\partial x$, and from
(\ref{eq:vel}) and (\ref{eq:poisson}) we have $-\laplacian u
=\partial\kern.3pt\omega/\partial y$, implying 
$\partial \omega/\partial y = -\partial p/\partial x$.
Similarly, $-\laplacian v=-\partial p/\partial y$ and $-\laplacian v =\partial
\omega/\partial x$, giving
$\partial\kern.3pt\omega/\partial x =\partial p/\partial y$.
These are the Cauchy--Riemann equations for $\omega$ and $p$ as functions of $x$
and $y$, implying that $\omega$ and $p$ are harmonic conjugates.  Therefore, the
system of equations~(\ref{eq:incomp}) is equivalent to the stream
function-vorticity formulation
\begin{equation}
-\laplacian\psi = \omega, \quad
-\laplacian\kern .8pt \omega = 0.\label{eq:laplace}
\end{equation}
Combining these equations shows that $\laplacian^2
\psi=0$, i.e., $\psi$ is biharmonic.
Note that the pressure has been eliminated from the problem, reducing
the original system of three equations in three unknowns to the biharmonic
equation in a single variable. By Theorem \ref{th:goursat}, one can write
\begin{equation}
\psi = \frac{1}{2\kern .3pt i} \bigl(\kern .7pt g(z)
- \overline{g(z)} + \overline{z}f(z) - z\overline{f(z)}\kern 1.6pt \bigr)
\label{eq:psieq}
\end{equation}
for some analytic functions $f(z)$ and $g(z)$, whereupon the velocity components become
\begin{equation}
u-iv = g'(z) + \overline{z} f'(z) - \overline{f(z)}.
\label{eq:velocities}
\end{equation}
Since $\psi$ is invariant under the transformations of Theorem
\ref{th:goursat}, so are $u$ and $v$. The vorticity is given by
\begin{displaymath}
-\omega = \laplacian \psi =
4\frac{\partial^2\psi}{\partial \overline{z}\kern 1pt\partial z} 
   = -2\kern .5pt i
\bigl( \kern .5pt f'(z) - \overline{f'(z)}\kern 1.6pt \bigr) =
4\kern .5pt\Im\kern -.8pt \left\{f'(z)\right\},
\end{displaymath}
and from the second equation of (\ref{eq:laplace}) we deduce that
\begin{equation}
p-i\kern .4pt\omega = 4 f'(z).
\end{equation}
Note that since the pressure is a real quantity defined up to a constant,
this equation indicates that $f'(z)$ is defined up to a real constant---the
number $\gamma$ in Theorem~\ref{th:goursat}.

In applications, eq.~(\ref{eq:velocities}) is the key to enforcing the
boundary conditions (\ref{eq:bcsa})--(\ref{eq:bcsb}).
Writing the velocities out individually gives
\begin{equation}
\label{eq:explicit}
u = \Re (\kern .8pt g' - f + \overline{z} f') , \quad
v = \Im (- g' - f - \overline{z} f') ,
\end{equation}
and it is by imposing these conditions, together with the 
condition (\ref{eq:psieq}) on boundary segments where
$\psi$ is known, that we will determine $f$ and $g$.

For our numerical examples, a phenomenon of particular
interest will be the appearance of {\em Moffatt eddies}
near corners where two solid boundaries meet with
straight sides at a fixed angle, conventionally
denoted by $2\alpha$~\cite[chap.~12]{daup,hasi,langdev}.
Building on earlier work by Dean and Montagnon~\cite{dean49}, 
Moffatt showed in 1964 that if $2\alpha$ is less than a critical value of
about $146^\circ$, then in principle one can expect a self-similar infinite
series of counter-rotating eddies to appear with
rapidly diminishing amplitude.  Physically, such eddies involve fluid
that is trapped near the corners, and the curves separating one eddy
from the next represent separation of the flow at the boundary, an
effect that would be absent if there were no viscosity (i.e., potential flow).
In polar coordinates $r$, $\theta$
Moffatt's asymptotic solution takes the form
\begin{equation}\label{eq:asymp}
\psi_\lambda(r,\theta) \sim
\kern .8pt \Re\bigl\{ A_0 \kern .8pt r^\lambda
   \bigl[\kern .7pt \cos{(\lambda\alpha)}\cos{((\lambda-2)\theta)} +
   \cos{((\lambda-2)\alpha)}\cos{(\lambda\theta)}\bigr]
\bigr\},
\end{equation}
where, to satisfy the no-slip conditions
$\psi_\lambda=\boldsymbol{\hat n}\cdot\gradient\psi_\lambda=0$ at $\theta=\pm\kern .3pt \alpha$,
the eigenvalue $\lambda$ must be a root of
\begin{equation}\label{eq:roots}
\sin{\left(2\alpha(\lambda-1)\right)} + (\lambda-1)\sin{(2\alpha)}=0.
\end{equation}
For cases with Moffatt eddies, $\lambda$ will be complex.  With
angle $2\alpha = 90^\circ$, each eddy is about 16 times smaller in
spatial scale than the last, with velocity amplitude about $2200$
times smaller, as we shall see in Figure~\ref{fig:ldc}.  The ratio
of stream function amplitudes is about 36,000.  As $\alpha\to 0$,
the ratio of space scales approaches~$1$, as we shall begin to see
in Figure~\ref{fig:wedge}, but the ratios of velocity or stream
function amplitudes never fall below about $350$.  (The limit of
an infinite half-strip is the so-called Papkovich--Fadle problem.)
In a laboratory, one would not expect to observe more than one
or two eddies, and they are also hard to resolve computationally,
especially by finite element methods unless careful attention is
given to mesh refinement~\cite{beckermao,cpr,kondrat,ronquist}.

\section{Numerical Method\label{sec:method}}

We now describe our method for solving
(\ref{eq:bih})--(\ref{eq:bcsb}) when $\Omega$ is a polygon
or curved polygon with $K$ corners $\{w_k\}_{k=1}^K$.
Following~\cite{gopal19}, the method consists of the rational
approximation of the analytic functions $f(z)$ and $g(z)$
of~(\ref{eq:goursat}) satisfying boundary conditions expressed
in terms of (\ref{eq:explicit}).  We choose $N$ basis functions
$\{\phi_j(z)\}_{j=1}^N$ with fixed poles, and we determine
the coefficients by solving a linear least-squares problem
in $M\gg 2N$ sample boundary points $\{z_i\}_{i=1}^M \subset
\partial\kern.5pt\Omega$.  Sometimes for simplicity we specify
parameters a priori, and in other cases, we increase $N$ adaptively
until the residual is below a given tolerance or the refinement
does not bring any improvement. Our implementation of the adaptive
version of the method is modeled on the publicly available code
$\texttt{laplace.m}$ \cite{tref20}, currently in version 6.
The algorithm is summarized in the listing labeled Algorithm 1.

\begin{algorithm}
\DontPrintSemicolon
\openup 1.5pt
\SetAlgoLined
\vrule width 0pt height 11pt
Define the boundary $\partial\kern.5pt \Omega$
and corners $w_1,\ldots,w_K$,\\[-1pt]
\indent\kern .3in the boundary data $h$, $k$, $\mathbf a$, $\mathbf g$, and 
the tolerance $\varepsilon$.\\
   \Repeat{\rm $\|Ax-b\|_2<\varepsilon$ or $N$ is too large or the error is growing.}{
      \vrule width 0pt height 10pt Define the basis of $N$ rational functions,\\
      Choose $M\gg 2N$ sample points $z_1,\dots,z_M$, clustered near the corners,\\
      Form the $2M\times 4N$ real matrix $A$ and right-hand side vector $b$,\\
      Solve the least-squares problem $Ax\approx b$ for $x$, after row weighting,\\
      \vrule width 0pt depth 4pt
      Compute residuals near each corner and increment the numbers of poles\\
   }
   \vrule width 0pt height 9pt
   Confirm accuracy by checking the error on a finer boundary mesh.\\
   \vrule width 0pt depth 6pt
   Construct functions to evaluate $u(z)$ and $\psi(z)$ at arbitrary
   points $z\in\Omega$.\\
   \caption{Lightning Stokes solver (adaptive version)\vrule width 0pt height 9.5pt depth 3pt}
\end{algorithm}

\subsection{Rational functions}
The rational approximations take the form
\begin{equation}
f(z) = \sum_{j=1}^N f_j \phi_j(z), \quad g(z) = \sum_{j=1}^N g_j \phi_j(z),
\end{equation}
where $f_j, g_j\in\complexes$ are complex coefficients
and the basis functions $\phi_j$ can be written as
\begin{equation}
\left\{\phi_j(z)\right\}_{j=1}^{N} = 
   \bigcup_{k=1}^K\left\{\frac{1}{z-\beta_{kn}}\right\}_{n=1}^{N_k} 
   \cup \left\{\kern .8pt p_n(z)\right\}_{n=0}^{N_0}.
\label{eq:allpoles}
\end{equation}
The functions in the left half of this union consist of simple poles
at preassigned locations clustered near the corners of $\Omega$, and
those on the right define a polynomial to handle ``the smooth part
of the problem.''  In~\cite{gopal19} these are called {\em Newman}
and {\em Runge} terms, respectively.  There are $4N$ real degrees of
freedom all together.

For the simple poles, we choose points $\beta_{kn}$ with tapered exponential
clustering near the corners $w_k$, following the formula introduced
in~\cite{gopal19} and analyzed in~\cite{clustering},
\begin{equation} \label{eq:poles}
\beta_{kn} = w_k + L\kern .5pt e^{i\theta_k} e^{-\sigma (\sqrt{N_k}-\sqrt{n}\kern .8pt )}, 
   \quad k=1,\ldots,K,\quad n=1,\ldots,N_k,
\end{equation}
with $\sigma=4$.  Here
$\theta_k$ is the angle of the
exterior bisector at the corner $w_k$ with respect
to the real axis and $L$ is a characteristic length scale
associated with $\Omega$.  For low-accuracy computations,
the functions $p_n$ can simply be taken as the monomials
$1, z, \dots, z^{N_0},$ but for higher accuracies it is
necessary to use discrete orthogonal polynomials
generated by the Vandermonde with Arnoldi process
described in~\cite{brubeck19} (Stieltjes orthogonalization), spanning the
same spaces but in a well-conditioned manner.
They are related by the $(n+2)$-term recurrence relation
\begin{equation} \label{eq:arnoldi}
h_{n+1,n} p_{n+1}(z) = z\kern .3pt p_n(z) - \sum_{k=0}^n h_{k,n} p_k(z),
   \quad n=0,\ldots,N_0-1,
\end{equation}
with $p_0(z)=1$,
where the numbers $h_{ij}$ are the entries of $H$, the upper-Hessenberg
matrix generated in the Arnoldi process.  

Unlike Laplace problems, Stokes problems invariably require
derivatives of the basis functions $\phi_j$, since $f'$
and $g'$ appear in (\ref{eq:explicit}).
Differentiating (\ref{eq:allpoles}) gives
\begin{equation}
\left\{\phi_j'(z)\right\}_{j=1}^{N} = 
   \bigcup_{k=1}^K\left\{\frac{-1}{(z-\beta_{kn})^2}\right\}_{n=1}^{N_k} 
   \cup \left\{\kern .8pt p_n'(z)\right\}_{n=0}^{N_0}.
\label{eq:allpolesderiv}
\end{equation}
If $p_n(z) = z^n$, then $p_n'(z) = n z^{n-1}$, but when Vandermonde
with Arnoldi orthogonalization is in play it is necessary to 
calculate $p_n'$ from the coefficients $h_{ij}$ of (\ref{eq:arnoldi}),
\begin{equation} \label{eq:arnoldidiff}
h_{n+1,n} p_{n+1}'(z) = p_n(z) + z\kern .3pt p_n'(z) - \sum_{k=0}^n h_{k,n} p_k'(z)
   \quad n=0,\ldots,N_0-1.
\end{equation}

In the description above, Arnoldi orthogonalization is applied to
the polynomial part of (\ref{eq:allpoles}) but not to the rational part, and experience
suggests that this often suffices for accurate computations.
However, our usual practice is to orthogonalize the
rational part of the approximation too.  This idea originates
with Yuji Nakatsukasa (unpublished) and is realized in the codes
{\tt VAorthog} and {\tt VAeval} listed in
the appendix (Figure~\ref{fig:arnoldicodes}).

\subsection{Linear least-squares system}

The boundary conditions \eqref{eq:bcsa}--(\ref{eq:bcsb})
are enforced at a set of $M \gg 2N$
sample points $\{z_i \}_{i=1}^M$, which are clustered near the corners
similarly as $\beta_{kn}$.  For simple computations we typically
cluster sample points by prescribing, say, 
{\tt x = tanh(linspace(-14,14,300))} for 300 clustered points
in the interior of the boundary segment $[-1,1]$.  For more careful adaptive
computations, sample points are introduced in step with the introduction
of clustered poles.  The boundary conditions are then imposed
by solving the real least-squares problem
\begin{equation} \label{eq:LS}
A x\approx b,
\end{equation}
with $A\in\reals^{2M\times 4N}$, $x\in\reals^{4N}$, and
$b\in\reals^{2M}$.
To describe the entries of $A$, it
is convenient to partition the system into $2\times 4$ blocks
$A_{ij}$ associated with sample points $z_i\in\partial\kern.5pt\Omega$ 
and basis functions $\phi_j$:
\begin{equation}
\begin{bmatrix}
A_{11} & \cdots & A_{1N} \\
\vdots & \ddots & \vdots \\
A_{M1} & \cdots & A_{MN}
\end{bmatrix}
\begin{bmatrix}
x_1 \\ \vdots \\ x_N
\end{bmatrix}
\approx
\begin{bmatrix}
b_1 \\ \vdots \\ b_M
\end{bmatrix}.
\label{eq:leastsq}
\end{equation}
The four columns of 
$A_{ij}\in\reals^{2\times 4}$ correspond to the real and imaginary
parts of the complex coefficients of $\phi_j$ in $f$ and $g$
associated with the $4\times 1$ block of unknowns
$x_j = (\Re f_j, \Im f_j, \Re g_j, \Im g_j)^T$,
and the two rows correspond to the boundary conditions applied at $z=z_i$
associated with the $2\times 1$ block of boundary values
$b_i=(U(z_i), V(z_i))^T$.  By (\ref{eq:explicit}), we have
\begin{displaymath}
\renewcommand*{\arraystretch}{1.4}
A_{ij} =
\begin{bmatrix}
\hphantom{-}\Re\{\overline{z_i}\kern .7pt \phi'_j(z_i)-\phi_j(z_i)\}
& \kern -4pt-\Im\{\overline{z_i}\kern .7pt \phi'_j(z_i)-\phi_j(z_i)\}
& \kern -2pt\hphantom{-}\Re\kern .8pt \phi'_j(z_i) 
& \kern -4pt -\Im\kern .8pt  \phi'_j(z_i) \\
-\Im\{\overline{z_i}\kern .7pt \phi'_j(z_i)+\phi_j(z_i)\} 
& \kern -4pt -\Re\{\overline{z_i}\kern .7pt \phi'_j(z_i)+\phi_j(z_i)\}
& \kern -2pt -\Im \kern .8pt \phi'_j(z_i)
& \kern -4pt -\kern .8pt \Re\kern .8pt \phi'_j(z_i)
\end{bmatrix}.
\end{displaymath}
A no-slip condition on a stationary boundary segment
corresponds simply to $b_i=(0,0)^T$.
If the boundary segment is moving tangentially at
speed $c$, as in a driven cavity,
one wants to set the normal and tangential velocities to $0$
and $c$, respectively.  (Alternatively, when the value of $\psi$
on the boundary is known a priori, it is better to impose this directly
as a Dirichlet boundary condition rather than setting the normal
velocity to~$0$.)  If the moving boundary segment is horizontal,
these velocity components are $\pm v$ and $\pm u$, and if it is
vertical, they are $\pm u$ and $\pm v$.  For boundary segments
oriented at other angles, some trigonometry is involved that
can conveniently be implemented by left-multiplying the matrices
$A_{ij}$ and vectors $b_i$ by a $2\times 2$ rotation matrix to
align the boundary with an axis.

\subsection{Row weighting}
In our MATLAB implementations, we solve \eqref{eq:LS} via the
backslash command.  First, we usually introduce row weighting,
for two reasons.  One is that our
grids are exponentially clustered near the corners, giving
these regions greater weight; we compensate for this effect
by multiplying row $i$ of (\ref{eq:LS}) by $|z_i-w_k|^{1/2}$,
where $w_k$ is the nearest corner of $\partial\kern .5pt \Omega$
to $z_i$.  In addition, $f'$ and $g'$ are usually discontinuous
at certain corners, so that pointwise convergence of the numerical
approximations is not an option.  Multiplying row $i$ by a second
factor of $|z_i-w_k|^{1/2}$ compensates for this second effect.
Combining the two scalings leads to our standard practice of
multiplying every row of $A$ and the right-hand side by the
appropriate factor $|z_i-w_k|$.

\subsection{Elimination of four redundant degrees of freedom}
According to Theorem~\ref{th:goursat}, $f$ and $g$ are uniquely
defined up to four real parameters,  giving the matrix $A$ as
defined so far a rank-deficiency of 4.\ \ Experiments show that
if one ignores this, there is often good convergence of $\psi$
anyway, even though $f$ and $g$ do not converge individually.
However, we usually eliminate these degrees of freedom by
appending four further rows to $A$ and to the right-hand side.
For the square lid-driven cavity of the next section, our normalization
conditions are arbitrarily chosen as
$f(0) = \Re \kern .8pt g(0) = \Re f(1) = 0$.

\subsection{Adaptivity}
When computing in adaptive mode,
we use a similar adaptive strategy to that of~\cite{gopal19},
which consists of solving a sequence of linear
least-squares problems of increasing size.
In brief, we associate each entry of
the residual vector $r=b-Ax$ with the nearest corner $w_k$ by defining the
index sets
\begin{displaymath}
\mathcal{I}_k = \Big\{1\leq i\leq M : k=\arg\min_{1\leq l \leq K} |w_l-z_i|\Big\},\quad k=1,\ldots, K.
\end{displaymath}
After solving one least-squares problem, we compute the local errors
\begin{displaymath}
R_k^2 = \sum_{j\in \mathcal{I}_k} r_j^2,\quad k=1,\ldots, K.
\end{displaymath}
For the next least-squares problem, we increase the number of poles at
each corner, with a greater increase near corners with larger errors.
The computation is terminated if the error estimated on a finer
boundary mesh is sufficiently small or $N$ is too large or the error is
not decreasing.

\section{Examples:~bounded cavities\label{sec:results}}
\begin{figure}
\vskip 5pt
\centering
\includegraphics[scale=.9]{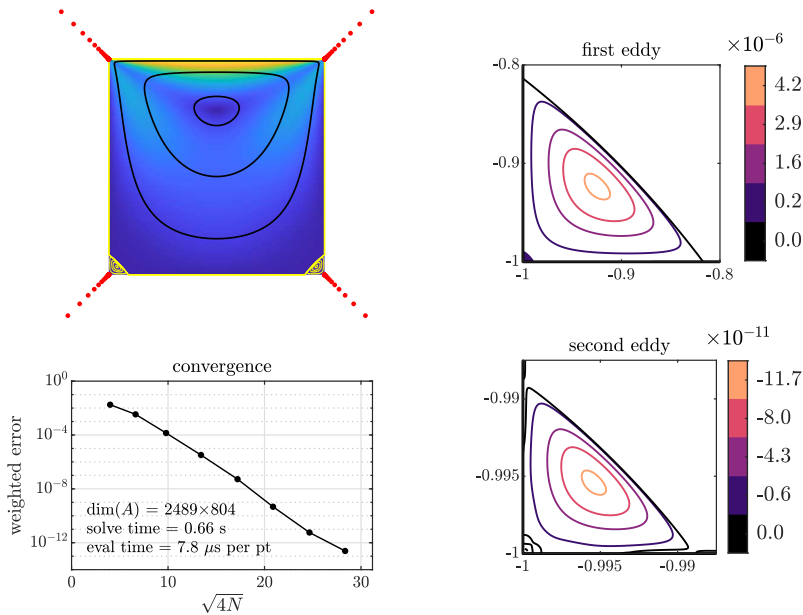}
\caption{Stokes flow in a square lid-driven cavity computed to
$12$-digit accuracy
in half a second on a laptop.  Stream contours and velocity magnitude (top
left), convergence (bottom left),
and close-ups of the first two Moffatt eddies in the lower-left corner (right).
The red dots mark the exponentially clustered
poles $(\ref{eq:allpoles})$ of the lightning approximation, $49$ at each upper
corner and\/ $39$ at each lower corner (with\/ $4$ poles outside the
plotting axes at each corner).}
\label{fig:ldc}
\end{figure} 

In this section and section~\ref{sec:unbounded}
we illustrate the lightning Stokes solver
applied to various model problems.  Each example
presents a contour plot of the stream function $\psi$ superimposed
on a color plot of the velocity magnitude $\sqrt{u^2+v^2}$, using
black contours for the main flow and yellow for the Moffatt eddies.
The poles $\beta_{kn}$ are marked by red dots. As shown by
plots in Examples~1 and~3, all the computations exhibit
root-exponential convergence
as a function of $4N$, the number of real degrees of
freedom, as expected from the theory of~\cite{gopal19}.  
The error measure in the plots is the maximum of the pointwise deviation from
the boundary conditions weighted by distance to the nearest corner.

In standard floating point arithmetic, we can often get
10--12 digits of accuracy in $\psi$ with 20--40 poles clustered
at each corner.  In this regime the floating point limits
are beginning to matter, and we have hand-tuned a
few parameters to get the desired accuracy.  Because of the
root-exponential convergence, on the other hand, one can
achieve 5--6 digits of accuracy with only 5--10
poles at each corner, with a reduction in linear algebra costs
in principle by a factor of $1/16$ or better depending on details
of boundary sampling and adaptivity.  (In practice it is usually
less because the computations are not in the asymptotic regime.)
These lower-accuracy computations are easy, typically requiring
no adaptivity or careful tuning of parameters.

\smallskip 
{\bf Example 1.  Lid-driven square cavity.} 
Figure~\ref{fig:ldc} presents a square
lid-driven cavity, where the fluid is enclosed between solid
walls one of which moves tangentially at speed~1~\cite{review}.
The domain is $\Omega=[-1,1]^2$, with boundary conditions $\psi
= 0$, $\psi_n = 1$ on the top and $\psi=\psi_n=0$ on the other
three sides.  Since $\psi_n$ is discontinuous at the top corners,
$\psi$ has singularities there.  At the bottom corners there are
weaker singularities associated with the formation of Moffatt
eddies, but $\psi_n$ is continuous, and our adaptive strategy puts
fewer poles at these corners.

The plots on the right in Figure~\ref{fig:ldc}
show close-ups of two eddies from the lower-left corner
of the flow domain, with contours in each case at levels 0, 0.05,
0.35, 0.65, and 0.95 relative to the maximum value of $|\psi|$.
Clearly the lightning method does a successful job of resolving
two eddies, even though the stream function in the second one has
amplitude of order $10^{-10}$.  The eddies have the same shapes
(on different space scales), opposite signs, and enormously
different amplitudes.

\begin{figure}
\vskip 17pt
\centering
\includegraphics[scale=.85]{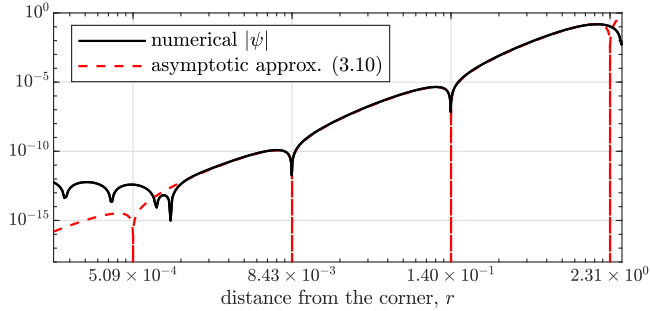}
\caption{Log-log plot of the stream function magnitude
$|\psi|$ along the diagonal of the lid-driven cavity.
The computed stream function matches the asymptotic approximation
$(\ref{eq:asymp})$ over an amplitude range of ten orders of magnitude.  The
constant $A_0$ of\/ $(\ref{eq:asymp})$ was determined by
least-squares fitting.}
\label{fig:ldc_loglog}
\end{figure}

Note that the Moffatt eddies have the same forms in the lower-left
and lower-right corners.  In this as in all Stokes flows, reversing
the signs of the boundary data just reverses the direction of the
flow, since \eqref{eq:bih} is linear.   In other words it doesn't
matter whether the lid is moving from left to right or right
to left.  This reversibility does not hold in Navier--Stokes
flows, with nonzero Reynolds number and a nonlinear governing
partial differential equation.\footnote{Reversibility of Stokes flows
is the reason microorganisms use flagella rather than fins to get
around.  See~\cite{lauga} for a review of this beautiful subject.}

Figure~\ref{fig:ldc_loglog} illustrates the power of asymptotic analysis,
and the accuracy of the lightning solver, with a plot of $|\psi|$ along
the $45^\circ$ line from the lower-left corner.  The numerical results,
displayed by the solid line, are compared with an asymptotic estimate
based on (\ref{eq:asymp}).  To obtain this estimate, we first calculated
the complex eigenvalue $\lambda \approx 3.740 + 1.119\kern .3pt i$ from
\eqref{eq:roots} with $2\alpha=90^\circ$, and the complex constant $A_0$
of (\ref{eq:asymp}) was then determined by least-squares fitting of the
$|\psi|$ data over the interval from the corner to the boundary between
the eddies and the main flow.  The curves show agreement of numerics
with the asymptotic theory to about 12 digits.

These lightning calculations can be quite simple.
The appendix illustrates this with a MATLAB code
for the solution of the lid-driven cavity problem with 24 poles fixed
near each corner and a polynomial of degree
$24$ (Figure~\ref{fig:code}); this runs in
0.2 s on the same laptop and computes $\psi(0) \approx 0.11790231118443$
with an error of $1.7\times 10^{-13}$.  If 24 is reduced to~6,
the time is 0.04~s and the error is $4.8\times 10^{-7}$.  For this
geometry, the polynomial term is not very important
to the solution: if its degree is reduced to~0, the errors
increase only slightly to $6.2\times 10^{-12}$ and
$7.4\times 10^{-7}$.

\begin{figure}
\vskip 12pt
\centering
\includegraphics[scale=.67]{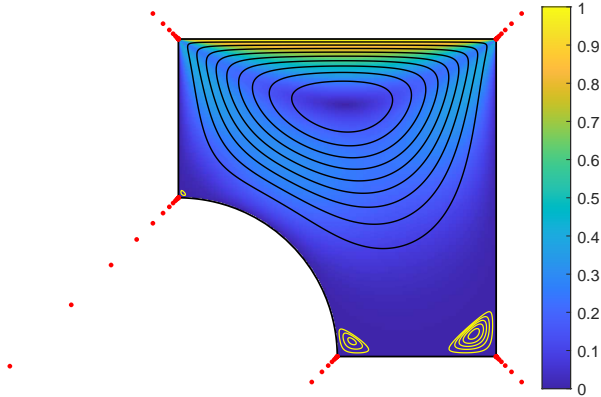}~~~~~~~~~~~
\caption{Lid-driven cavity with a quarter-circular exclusion,
with $\psi$ calculated to $8$ digits in $0.3$ seconds.}
\label{fig:ldccorner}
\end{figure} 

\smallskip 
{\bf Example 2.  Lid-driven cavity with circular exclusion.}  A second
lid-driven cavity is shown in Figure~\ref{fig:ldccorner}, now with a
quarter-circular arc removed from the boundary.
The computation is similar, and
with polynomial degree $40$ and $20$ poles at each of the three
corners, we get the value $\psi(0)\approx -0.0599323802$ to
accuracy $6\times 10^{-9}$ in $0.3$ s.

\begin{figure}
\indent\kern 1.25in\includegraphics[scale=.9]{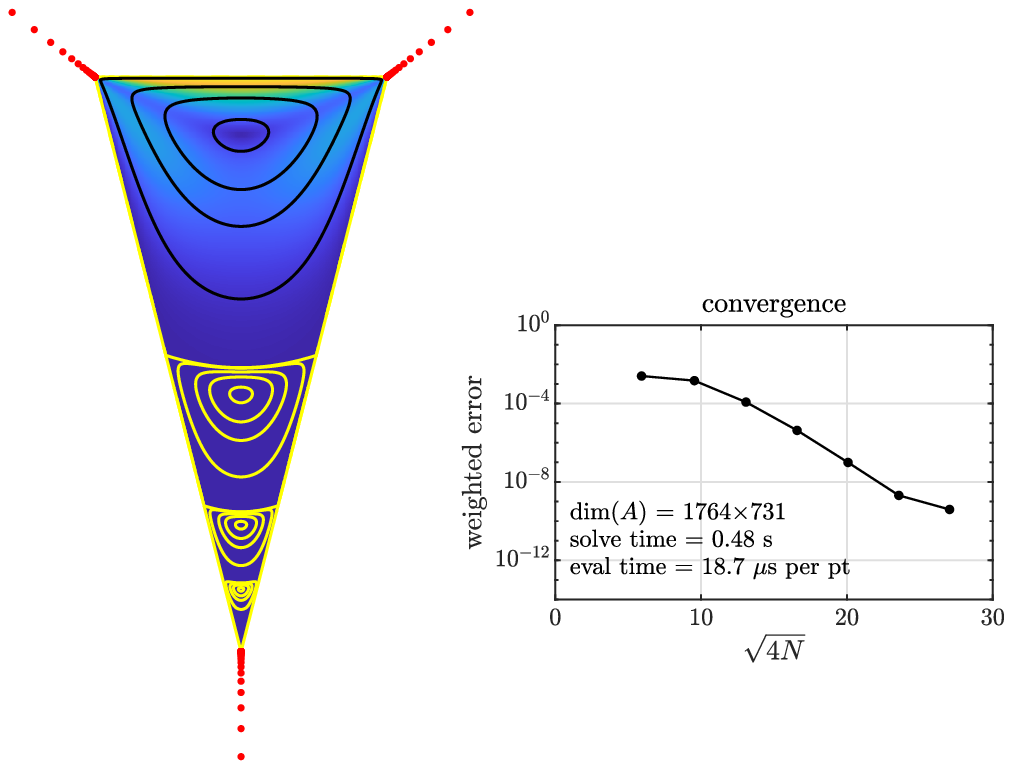}
\pput(-339,32){\includegraphics[scale=.56]{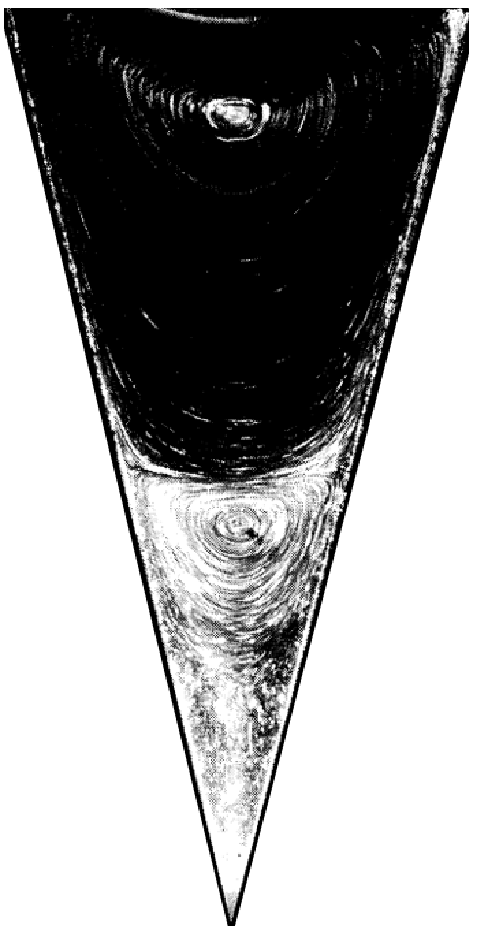}}
\caption{Stokes flow in a triangular lid-driven cavity 
of vertex angle $2\alpha=28.5^\circ$.  There are $49$ poles at
each corner, of which\/ $7$ at the lower corner and\/ $5$ at each upper
corner lie outside the plotting axes.  The computed result matches
Taneda's experiment from $1979$~{\rm \cite{taneda79}}.}
\label{fig:wedge}
\end{figure}

\smallskip
{\bf Example 3.  Lid-driven triangular cavity.} Our
third example, shown in Figure~\ref{fig:wedge}, is an isoceles
triangle with unit side length and vertex angle $2\alpha = 28.5^\circ$,
corresponding to the complex eigenvalue $\lambda = 9.485 + 4.434i$ in
(\ref{eq:roots}).  The stream function is resolved to about 10 digits in
half a second.  The sharper corner makes the eddies decrease
more slowly than before, with the ratio of successive amplitudes $|\psi|$
now about $830$.  Three eddies can be observed without
the need of a close-up.

This particular setup was chosen for comparison with experimental images
from Taneda~\cite{taneda79}, which are well known from Van Dyke's {\em
Album of Fluid Motion}~\cite{vandyke82}.   In his experiment, Taneda drove
a flow of silicone oil at Reynolds number $0.17$ by a rotating cylinder,
with flow visualization by aluminum powder and a photographic exposure
time of 90 minutes.  Only one Moffatt eddy is clearly visible in the
experiment, as the velocities are nearly zero close to the corner.
According to Taneda, ``the reason is that the relative intensity of
successive vortices is the order $10^3$, and therefore the photographic
exposure necessary to visualize two successive vortices is about $10^3$
times that necessary to visualize a single vortex.''  As far back
as 1988, R\o nquist successfully computed Taneda's flow to a similar
accuracy as ours using spectral element methods:~30 elements, each of
order 8~\cite[Fig.~12]{ronquist}.

A comparison of the numerical results with the asymptotic
approximation (\ref{eq:asymp}) as in Figure~\ref{fig:ldc_loglog}
(not shown) again indicates close agreement down to
$|\psi|\approx 10^{-12}$.

\section{Behavior outside the problem domain}\label{sec:analytic}Finite difference
and finite element methods represent the solution to a PDE directly in
the problem domain $\Omega$, and integral equations methods represent
it via density functions supported on the boundary $\partial\kern
.5pt \Omega$.  In neither case is anything computed outside $\Omega$.
Rational approximation methods, however, represent the solution by a
function that is defined throughout the plane.  In regions of analyticity,
the computed approximation will typically provide an accurate analytic
continuation across the boundary, whereas near branch points, it will
feature a string of poles outside $\Omega$ that approximate
a branch cut.  For discussions of these effects
see the conclusions section of~\cite{gopal19} and Figure~3.2
of~\cite{loglight}.

\begin{figure}
\vskip 35pt
\begin{center}
\kern -12pt\includegraphics[scale=.9]{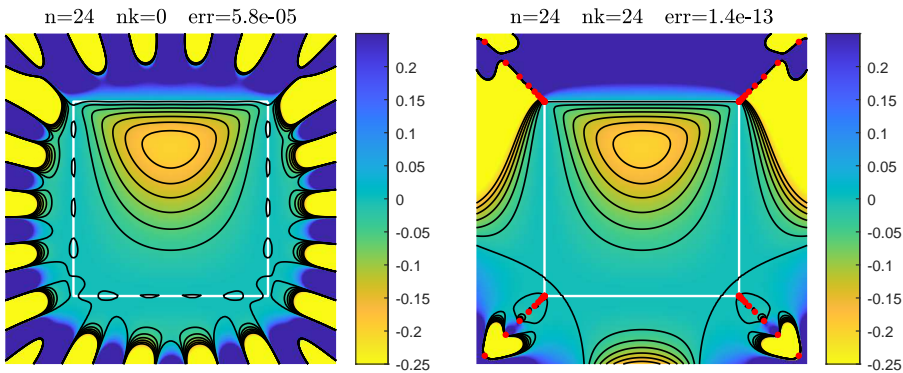}
\end{center}
\caption{\label{fig:analytic}Approximations to the real function
$\psi(z)$ for the lid-driven unit square cavity shown in the
larger square $[-1.7,1.7]^2$.  On the left,
a degree $24$ polynomial approximation, and on the right, a rational
approximation with, in addition, $24$ poles near each corner.
The strings of poles approximate branch cuts needed to approximate
the Goursat functions $f$ and $g$.
The cavity boundary, marked in white, does not show
up as a zero contour on the sides and the bottom
because $\psi$ has a local minimum or maximum there rather
than a change of sign.}
\end{figure}

For illustration, consider the biharmonic equation
in the lid-driven square cavity.  First,
Figure~\ref{fig:analytic} plots the computed approximation to
the real function $\psi(z)$, not just on $\Omega$ itself but on
the larger square $[-1.7,1.7]^2$.  (In contrast to our other plots, 
here the colors represent $\psi$, not $\sqrt{u^2 + v^2}$.)
he image on the left shows a
degree $24$ polynomial approximation, with the familiar circus
tent effect at larger values of $|z|$.  Such approximations
are good for getting a few digits of accuracy, at least away from
singular corners, and the error figure in the title reports 4-digit
accuracy of the quantity $\psi(0) \approx -0.11790231118443$.
However, convergence will be very slow as the degree is increased.
The image on the right shows the result when $24$ poles are also
included clustered near each corner.  Now $\psi(0)$ is accurate
to 13 digits, and the contours give some insight into how this is
achieved.  Although the functions $f$ and $g$ from which $\psi$ is
derived have branch points at all four corners, good approximations
interior to $\Omega$ are found thanks to the clustered poles.
Note that the sides and the bottom of $\partial \kern .5pt
\Omega$, marked in white, do not show up as black contours at
level $\psi=0$.  This is because $\psi$ has a local maximum rather
than a simple zero there (or a minimum, inside the Moffatt eddy),
which goes undetected by MATLAB's contour finding algorithm.
By contrast the separation lines between the main flow and the
Moffatt eddies represent sign changes and do appear in black.

\begin{figure}
\rlap{\includegraphics[scale=.96]{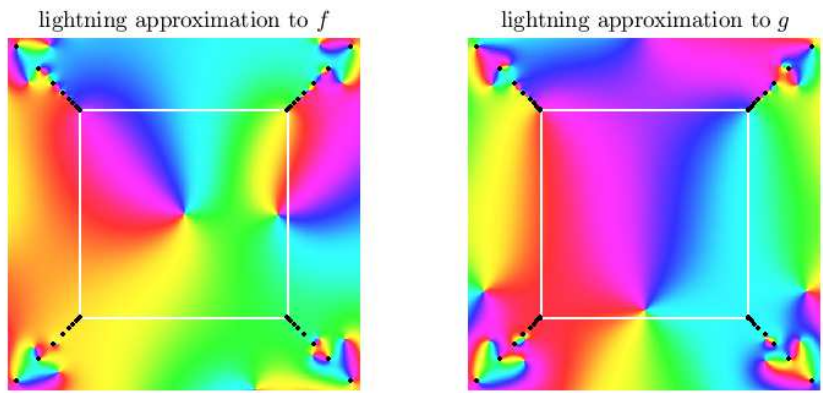}}
\vskip -142pt
\caption{\label{fig:analyticfg}Phase portraits of the computed
approximations to $f$ and $g$
in the same situation as the right image of Figure~$\ref{fig:analytic}$,
with $24$ clustered poles at each corner (now marked in black rather than
red).  The approximations to both
$f$ and $g$ show approximate branch cuts at the upper corners.
The singularities at the lower corners are weaker.}
\end{figure}

\begin{figure}
\rlap{\includegraphics[scale=.96]{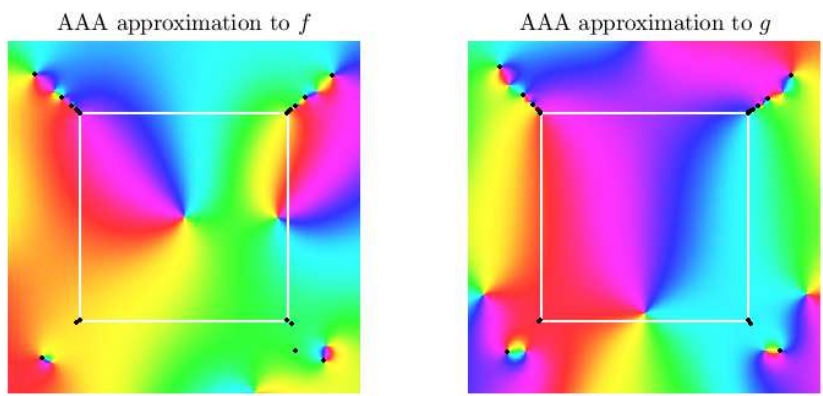}}
\vskip -142pt
\caption{\label{fig:analyticfgAAA}Phase
portraits of AAA approximations with tolerance $10^{-8}$ to the
functions of Figure~{\rm\ref{fig:analyticfg}}.  The poles are
determined adaptively by the AAA algorithm and appear clustered near
the corners, with $17$--\kern .8pt$18$ poles near each upper corner and
just $4$--\kern .8pt$5$ poles near each lower corner, where
the singularities are weaker.}
\end{figure}

To see more, we can examine the analytic functions $f$ and $g$
from which $\psi$ is constructed.  Figure~\ref{fig:analyticfg}
shows {\em phase portraits}~\cite{wegert}, with each
pixel colored according to the phase of the function value
there: red, yellow, green, cyan, blue, magenta for
phases $0,\pi/3,2\pi/3,\dots, 5\pi/3$.  Perhaps the first thing
one notices in these images is that $f$ has a zero at $z=0$.
This is not important, being just a result of our normalization.
The second feature one notices is the important one, that the
approximations to $f$ and $g$ both show what look like branch cuts
emanating from the upper corners.  This is a familiar effect in
rational approximation, and it illustrates why it is important
for poles to be placed near these corners.  One cannot say that
the exact functions $f$ and $g$ themselves have branch cuts, for
they are only defined on the unit square.  But clearly
they have branch {\em points} at $(1,1)$ and $(-1,1)$.
They also have branch points at the
lower corners, but these are weaker and not visible in the images.

In Figure~\ref{fig:analyticfg}, the poles are situated exactly
where we placed them according to (\ref{eq:poles}), with 24 poles
at each corner.  Now that we have good approximations to $f$ and
$g$ on $\partial\kern .5pt\Omega$, however (hence also on $\Omega$
by the maximum modulus principle), it is possible to re-approximate
these functions by new rational functions with free poles and see
what structures emerge.  Figure~\ref{fig:analyticfgAAA} shows what
happens when this is done by AAA approximation~\cite{AAA} with a
tolerance of $10^{-8}$.\footnote{AAA approximations were computed
in Chebfun with the {\tt aaa} command, and Chebfun {\tt phaseplot}
was used for the phase portraits.} Although the poles have not
been fixed in advance, we nevertheless find $17$ or $18$ poles
clustering exponentially near each upper corner, for both $f$
and $g$.  At the lower corners in this experiment, $f$ has 4 poles
and $g$ has 5.  These configurations confirm that the lightning
method does a good job of capturing the true structure of the
singularities of $f$ and $g$, and that there are singularities
at all four corners but the lower ones are weaker.

\section{Examples:~channel flows\label{sec:unbounded}}

We now show a pair of flows in channels as opposed to bounded
cavities.  Physically, the channels are infinite, and we approximate
them on bounded computational domains
by specifying the stream function upstream and downstream.  For the
best results here we found it effective to remove the row weighting
but scale all the columns of $A$ to have equal 2-norms.

\begin{figure}
\vskip 37pt
\begin{center}
\includegraphics[scale=1]{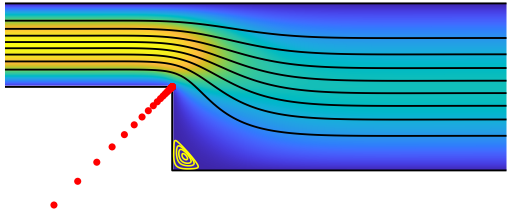}
\end{center}
\caption{Stokes flow over a step, with the infinite
channel approximated by a finite domain with explicit
inflow and outflow conditions.}
\label{fig:step}
\end{figure}

\smallskip
{\bf Example 4.  Square step.} Figure~\ref{fig:step} shows
the flow over a square step, with the reentrant corner at $(0,0)$ and
the salient corner at $(0,-1)$.  Our computational
domain introduces additional upstream corners at $(-2,0)$ and $(-2,1)$
and downstream corners at $(4,-1)$ and $(4,1)$.  (As always, thanks
to reversibility of Stokes flow, there is no true distinction
between upstream and downstream.)  Fully-developed parabolic
velocity profiles are imposed at the inflow and outflow such
that mass flow is conserved, $(u,v)=(1-(2y-1)^2,0)$ at $x=-2$
and $(u,v)=((1-y^2)/2,0)$ at $x=4$.  That is, the stream function
values are $\psi(-2,y) = 2y^2 - (4/3)y^3$ at the inflow and
$\psi(4,y) = 1/3 + y/2 - y^3/6$ at the outflow, and we complete
the boundary conditions by specifying $\psi_0=2/3$ on the top of
the channel, $\psi_0=0$ on the bottom, and $\psi_n= 0$ everywhere.
We normalize the Goursat functions by $f(0.5) = 0$,
$\Re \kern .7pt g(0.5) = 0$, and $\Re f(1) = 0$.

Figure~\ref{fig:step} was generated by a code with polynomial degree $40$
and $80$ poles clustered at the reentrant corner.
The computation runs in about 0.6~s on
our laptop and computes $\psi(z)$ to about 6 digits of accuracy:
for $z = 1,2,3$ the computed values are $0.259289$, $0.329814$, and $0.333990$.
(It surprised us that this last number could be greater than the centerline
value of $1/3$, but this seems to be genuine.  If the length of the
computational domain is increased to better approximate an
infinite channel, $\psi(3)$ increases slightly to $0.334054$.)
\begin{figure}
\indent\kern -19pt\rlap{\includegraphics[scale=.95]{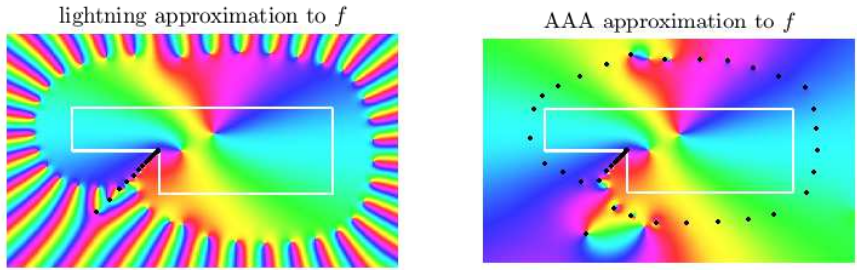}}
\vskip -192pt
\caption{Phase portraits of the Goursat function $f$ for the
flow over a step and of its AAA rational approximation with
tolerance $10^{-8}$.  The black dots marking poles confirm that
only the reentrant corner has a signiciant singularity.}
\label{fig:stepphase}
\end{figure} 

Following Figure~\ref{fig:analyticfgAAA}, Figure~\ref{fig:stepphase}
shows phase portraits of the lightning approximation to the Goursat
function $f$ for flow over the step, and its AAA approximant with
tolerance $10^{-8}$.

\begin{figure}
\vskip 55pt
\begin{center}
\kern -9pt\includegraphics[scale=1.1]{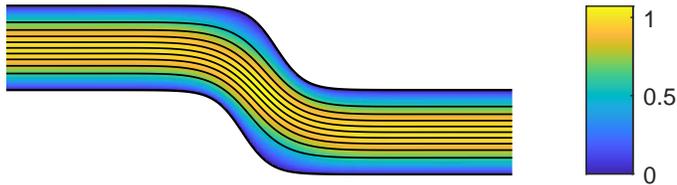}
\end{center}
\caption{Stokes flow in a smooth bent channel, computed
to $10$-digit accuracy by a polynomial expansion of degree $300$.}
\label{fig:bent}
\end{figure}

\smallskip
{\bf Example 5.  Bent channel.} Figure~\ref{fig:bent} shows
a flow in a smooth channel with a bend in it: the upper
boundary is $y = 1/2 - \tanh(3(x+.2))/2$ and the lower boundary is the
same curve rotated by $\pi$.
Since the walls are smooth, we use a purely
polynomial expansion of degree $n=300$.  The speed of flow is
$1$ in the middle of the channel at inflow and outflow, and 
at the center point $z = 0$ the computation gives speed $1.072187183679$,  
as compared with the true value of about 
$1.072187183704$.

\section{Discussion\label{sec:conclusion}}
The association of biharmonic functions and Stokes flow with
analytic functions goes back to the beginnings of these subjects,
and in various forms it has served as the basis of computations,
though perhaps not as frequently as one might have expected.
For example, a {\em stokeslet\/} is a biharmonic function with a
point singularity~\cite{daup,hasi,langdev}, and superpositions of stokeslets
have been used for numerical computations~\cite{cortez,fcb,lauga}.  What is
distinctive about the lightning approach to these problems
is the use of 
strings of singularities exponentially clustered near corners.
This leads to root-exponential convergence with respect to the
number of free parameters~\cite{gopal19,clustering}, an effect
first identified in rational approximation theory by Newman in
1964~\cite{newman64}.  With analysis of the singular behavior at
each corner, one could develop methods for certain problems with
even faster exponential convergence, but that requires case-by-case
analysis.  A more general approach to increasing convergence
rates, the ``log-lightning'' method based on reciprocal-log
approximations~\cite{loglight}, has not yet been investigated
for biharmonic problems.

The advantages of the lightning method are its simplicity,
its great speed and accuracy, and the fact that it delivers a
global representation of the solution (perfectly smooth
and exactly biharmonic).  A typical problem may be solved in less
than a second, with the solution in the form of a function handle
whose evaluation takes 5--50 microseconds per point.  We also
feel it is a worthwhile contribution conceptually, highlighting
as it does the analyticity properties that are such a fundamental
feature of these problems.

Following Goursat, our method explicitly constructs complex analytic
functions $f$ and $g$.  However, it has been known for a long
time that there are alternative formulations involving real
harmonic functions, associated with the names of Almansi and
Krakowski and Charnes:
here the representation
$\Im \kern -.8pt \left\{ \overline{z}f(z) + g(z)\right\}$
of (\ref{eq:goursat}) becomes $h(z)+xk(z)$, $h(z) +yk(z)$, or
$h(z) + r^2k(z)$, where $h$ and
$k$ are real and harmonic and $r = |z|$~\cite{fosdick,karfair88,saka}.
So far as we are aware, these representations could serve equally
well as the basis of a lightning solver.

We see two drawbacks to the lightning method.  One is that issues
of numerical stability are never far away and only partially
understood.  The Vandermonde with Arnoldi orthogonalization we
have described goes a good way to eliminating stability problems
in many cases, enabling computed results sometimes with close to
machine accuracy, but in other cases we find ourselves stuck at
10 digits or less.  We hope that further investigation will lead
to better undertanding and consequent improvements, as indeed
happened in the case of the original lightning Laplace method,
which was presented in~\cite{gopal19} before the introduction of
Vandermonde with Arnoldi~\cite{brubeck19}.

The other drawback is one that arises widely in computations
based on complex analysis, which is related to the ``crowding
phenomenon'' in numerical conformal
mapping~\cite[Thms.~2--5]{gtcm}.  Though the theorems may tell
us that exponentially clustered poles at each singularity ensure
root-exponential convergence, the constants associated with the
remaining, analytic part of the problem may be quite adverse for
domains involving certain distortions.  To combat such effects
one may employ representations of the Goursat functions $f$ and
$g$ that go beyond the basic prescription of clustered poles
+ polynomial, and this is a matter for further investigation.
It may prove fruitful to develop a variant of our algorithm that
combines exponentially clustered poles near singularities with
non-clustered poles spaced along curves, as in the more usual
method of fundamental solutions~\cite{karfair87,karfair88,saka}.
Indeed, the distribution of poles in the AAA approximation on
the right of Figure~\ref{fig:stepphase} seems to point in this direction.

The geometries we have considered here have been limited, but we
believe closely related methods will work for other cases too.
One prospect would be truly unbounded domains, such as infinite
channels, treated without artificial inflow and outflow conditions.
Complex variables and rational functions generally cope with
unbounded domains easily, so we are confident that this kind of
problem should be accessible.  Multiply-connected domains should
also be treatable without great changes.  At first this may seem
daunting, since the Goursat functions $f$ and $g$ will have to be
multi-valued~\cite{fosdick},
but in fact just a single complex logarithm will
be needed in each hole, as pointed out by Axler in his
wonderfully clear discussion of the ``logarithmic conjugation theorem''~\cite{axler}.
Series methods for Laplace problems with smooth holes are
considered in~\cite{series}, and holes with corners and lightning
terms appear in~\cite{conf}.   A third geometrical variation could
be problems with free boundaries~\cite[sec.~3.5]{ockendon95}.
And, of course, one may treat problems of elasticity as well as
fluid mechanics.

Extension of these methods to three dimensions is undoubtedly
possible, but not straightforward, and we have no view as to
whether something practical can be developed here or not.

For 2D Laplace problems, a lightning solver is available as a
MATLAB code that treats a wide variety of problems with a simple
interface~\cite{tref20}.  At present there is no comparable
software for the biharmonic case, so prospective users of this
method are recommended to begin by adapting the codes listed
in the appendix.

\section*{Appendix}
\begin{figure}
\label{fig:code}
{\footnotesize
\begin{verbatim}
      % Setup
        w = [1+1i; -1+1i; -1-1i; 1-1i];               % corners
        m = 300; s = tanh(linspace(-16,16,m));        % clustered pts in (-1,1)
        Z = [1i-s -1-1i*s -1i+s 1+1i*s].';            % boundary pts
        top = 1:m; lft = m+1:2*m; bot = 2*m+1:3*m;    % indices
        rgt = 3*m+1:4*m; sid = [lft rgt];
        n = 24; np = 24;                              % poly. deg.; # poles per corner
        dk = 1.5*cluster(np); d1 = 1+dk;              % tapered exponential clustering
        Pol = {w(1)*d1, w(2)*d1, w(3)*d1, w(4)*d1};   % the poles
        Hes = VAorthog(Z,n,Pol);                      % Arnoldi Hessenberg matrices
      
      % Boundary conditions
        [A1,rhs1,A2,rhs2,PSI,U,V] = makerows(Z,n,Hes,Pol);
        A1(top,:) = PSI(top,:); rhs1(top) = 0;        % Top:    psi = 0,
        A2(top,:) =   U(top,:); rhs2(top) = 1;        %         U = 1
        A1(bot,:) = PSI(bot,:); rhs1(bot) = 0;        % Bottom: psi = 0,
        A2(bot,:) =   U(bot,:); rhs2(bot) = 0;        %         U = 0
        A1(sid,:) = PSI(sid,:); rhs1(sid) = 0;        % Sides:  psi = 0,
        A2(sid,:) =   V(sid,:); rhs2(sid) = 0;        %         V = 0
        A = [A1; A2]; rhs = [rhs1; rhs2];             % put the halves together
      
      % Solution and plot
        [A,rhs] = rowweighting(A,rhs,Z,w);            % row weighting (optional)
        [A,rhs] = normalize(A,rhs,0,1,Hes,Pol);       % normalize f and g (optional)
        c = A\rhs;                                    % solve least-squares problem
        [psi,uv,p,omega,f,g] = makefuns(c,Hes,Pol);   % make function handles
        plotcontours(w,Z,psi,uv,Pol)                  % contour plot
        ss = 'Error in psi(0) = %6.1e\n';
        fprintf(ss,psi(0)+.117902311184435)           % accuracy check
\end{verbatim}
\par}
\caption{MATLAB program for the lid-driven square cavity, an example of a
ten digit algorithm as defined in~{\rm \cite{tda}}.  This code computes
the solution accurate to $13$ digits (as measured by $\psi(0)$) in $0.2$ secs on a
laptop and produces a plot in a further $3.5$ secs.  It calls the
functions {\tt VAorthog} and\/ {\tt VAeval} listed in Figure~{\rm\ref{fig:arnoldicodes}} for
Vandermonde with Arnoldi orthogonalization and the functions listed in
Figures~{\rm\ref{fig:morecodes}} and~{\rm\ref{fig:contours}}
for the remaining computation and plotting.}
\end{figure}
The method we have presented can be realized in
short MATLAB codes.  As a template, Figure~\ref{fig:code} presents an example
for the square lid-driven cavity, which readers may adapt to their
own problems.  This code calls a pair of functions for Vandermonde with Arnoldi
orthogonalization, shown in Figure~\ref{fig:arnoldicodes}, and additional functions for
setting up the algorithm and plotting, shown in Figures~\ref{fig:morecodes} and~\ref{fig:contours}.

\begin{figure}
\label{fig:arnoldicodes}
{\footnotesize
\begin{verbatim}
     function [Hes,R] = VAorthog(Z,n,varargin)  % Vand.+Arnoldi orthogonalization
     %  Input:   Z = column vector of sample points
     %           n = degree of polynomial (>= 0)
     %         Pol = cell array of vectors of poles (optional)
     % Output: Hes = cell array of Hessenberg matrices (length 1+length(Pol))
     %           R = matrix of basis vectors
     M = length(Z); Pol = []; if nargin == 3, Pol = varargin{1}; end
     % First orthogonalize the polynomial part
     Q = ones(M,1); H = zeros(n+1,n);
     for k = 1:n       
        q = Z.*Q(:,k);
        for j = 1:k, H(j,k) = Q(:,j)'*q/M; q = q - H(j,k)*Q(:,j); end 
        H(k+1,k) = norm(q)/sqrt(M); Q(:,k+1) = q/H(k+1,k);
     end
     Hes{1} = H; R = Q;
     % Next orthogonalize the pole parts, if any
     while ~isempty(Pol)
        pol = Pol{1}; Pol(1) = [];
        np = length(pol); H = zeros(np,np-1); Q = ones(M,1);
        for k = 1:np       
           q = Q(:,k)./(Z-pol(k));
           for j = 1:k, H(j,k) = Q(:,j)'*q/M; q = q - H(j,k)*Q(:,j); end 
           H(k+1,k) = norm(q)/sqrt(M); Q(:,k+1) = q/H(k+1,k);
        end
        Hes{length(Hes)+1} = H; R = [R Q(:,2:end)];
     end
\end{verbatim}
\vskip 9pt
\begin{verbatim}

     function [R0,R1] = VAeval(Z,Hes,varargin)  % Vand.+Arnoldi basis construction
     %  Input:   Z = column vector of sample points
     %         Hes = cell array of Hessenberg matrices
     %         Pol = cell array of vectors of poles, if any
     % Output:  R0 = matrix of basis vectors for functions
     %          R1 = matrix of basis vectors for derivatives
     M = length(Z); Pol = []; if nargin == 3, Pol = varargin{1}; end
     % First construct the polynomial part of the basis
     H = Hes{1}; Hes(1) = []; n = size(H,2); 
     Q = ones(M,1); D = zeros(M,1);
     for k = 1:n
        hkk = H(k+1,k);
        Q(:,k+1) = ( Z.*Q(:,k) - Q(:,1:k)*H(1:k,k)          )/hkk;
        D(:,k+1) = ( Z.*D(:,k) - D(:,1:k)*H(1:k,k) + Q(:,k) )/hkk;
     end
     R0 = Q; R1 = D;
     % Next construct the pole parts of the basis, if any
     while ~isempty(Pol)
        pol = Pol{1}; Pol(1) = [];
        H = Hes{1}; Hes(1) = []; np = length(pol); Q = ones(M,1); D = zeros(M,1);
        for k = 1:np
           Zpki = 1./(Z-pol(k)); hkk = H(k+1,k);
           Q(:,k+1) = ( Q(:,k).*Zpki - Q(:,1:k)*H(1:k,k)                   )/hkk;
           D(:,k+1) = ( D(:,k).*Zpki - D(:,1:k)*H(1:k,k) - Q(:,k).*Zpki.^2 )/hkk;
        end
        R0 = [R0 Q(:,2:end)]; R1 = [R1 D(:,2:end)];
     end
\end{verbatim}
\par}
\caption{MATLAB codes for Vandermonde with Arnoldi orthogonalization and
evaluation, adapted from~{\rm\cite{brubeck19}}.  The polynomial part
and each string of poles are orthogonalized separately,
with the resulting Hessenberg matrices of coefficients stored
in the cell array {\tt Hes}.}
\end{figure}

\begin{figure}
\label{fig:morecodes}
{\footnotesize
\begin{verbatim}

     function d = cluster(n)   % n points exponentially clustered in (0,1]
     nc = ceil(n); d = exp(4*(sqrt(nc:-1:1)-sqrt(nc))); 
\end{verbatim}
\vskip 9pt
\begin{verbatim}

     function [A1,rhs1,A2,rhs2,PSI,U,V] = makerows(Z,n,Hes,varargin)
     Pol = []; if nargin == 4, Pol = varargin{1}; end
     [R0,R1] = VAeval(Z,Hes,Pol);
     M = length(Z); N = 4*size(R0,2);
     cZ = spdiags(conj(Z),0,M,M);                    % conj(Z)
     PSI = [cZ*R0 R0]; PSI = [imag(PSI) real(PSI)];  % eq (2.1): stream function
     U = [cZ*R1-R0 R1]; U = [real(U) -imag(U)];      % eq (3.9): horizontal velocity
     V = [-cZ*R1-R0 -R1]; V = [imag(V) real(V)];     % eq (3.9): vertical velocity
     A1 = zeros(M,N); rhs1 = zeros(M,1);
     A2 = zeros(M,N); rhs2 = zeros(M,1);
\end{verbatim}
\vskip 9pt
\begin{verbatim}

     function [A,rhs] = normalize(A,rhs,a,b,Hes,varargin)
     % Impose conditions at z=a and z=b to make f and g unique
     Pol = []; if nargin == 6, Pol = varargin{1}; end
     rhs = [rhs; 0; 0; 0; 0];                       
     z = a; r0 = VAeval(z,Hes,Pol); zero = 0*r0;
     A = [A; real([r0 zero]) -imag([r0 zero])];      % Re(f(a)) = 0
     A = [A; imag([r0 zero])  real([r0 zero])];      % Im(f(a)) = 0
     A = [A; real([zero r0]) -imag([zero r0])];      % Re(g(a)) = 0
     z = b; r0 = VAeval(z,Hes,Pol);
     A = [A; real([r0 zero]) -imag([r0 zero])];      % Re(f(b)) = 0
\end{verbatim}
\vskip 9pt
\begin{verbatim}

     function [A,rhs] = rowweighting(A,rhs,Z,w)      % row weighting
     dZw = min(abs(Z-w.'),[],2);                     % distance to nearest corner
     wt = [dZw; dZw];
     M2 = 2*length(Z); W = spdiags(wt,0,M2,M2);
     A = W*A; rhs = W*rhs;
\end{verbatim}
\vskip 9pt
\begin{verbatim}

     function [psi,uv,p,omega,f,g] = makefuns(c,Hes,varargin)  % make function handles
     Pol = []; if nargin == 3, Pol = varargin{1}; end
     cc = c(1:end/2) + 1i*c(end/2+1:end);
     reshaper = @(str) @(z) reshape(fh(str,z(:),cc,Hes,Pol),size(z));
       psi = reshaper('psi');    uv = reshaper('uv');   p = reshaper('p');
     omega = reshaper('omega');   f = reshaper('f');    g = reshaper('g');

     function fh = fh(i,Z,cc,Hes,Pol)
     [R0,R1] = VAeval(Z,Hes,Pol);
     N = size(R0,2); 
     cf = cc(1:N); cg = cc(N+(1:N));
     switch i
        case   'f'  , fh = R0*cf;                                % eq (4.1): f
        case   'g'  , fh = R0*cg;                                % eq (4.1): g
        case  'psi' , fh = imag(conj(Z).*(R0*cf) + R0*cg);       % eq (2.1): psi
        case   'uv' , fh = Z.*conj(R1*cf) - R0*cf + conj(R1*cg); % eq (3.7): u+1i*v
        case   'p'  , fh = real(4*R1*cf);                        % eq (3.8): pressure
        case 'omega', fh = imag(-4*R1*cf);                       % eq (3.8): vorticity
     end
\end{verbatim}
\par}
\caption{Stokes discretization codes.}
\end{figure}

\begin{figure}
\label{fig:contours}
{\footnotesize
\begin{verbatim}

     function plotcontours(w,Z,psi,uv,varargin)   % contour plot
     Pol = []; if nargin == 5, Pol = varargin{1}; end
     MS = 'markersize'; LW = 'linewidth'; CO = 'color';
     x1 = min(real(Z)); x2 = max(real(Z)); xm = mean([x1 x2]); dx = diff([x1 x2]);
     y1 = min(imag(Z)); y2 = max(imag(Z)); ym = mean([y1 y2]); dy = diff([y1 y2]);
     dmax = max(dx,dy); nx = ceil(200*dx/dmax); ny = ceil(200*dy/dmax);
     x = linspace(x1,x2,nx); y = linspace(y1,y2,ny);
     [xx,yy] = meshgrid(x,y); zz = xx + 1i*yy;
     inpolygonc = @(z,w) inpolygon(real(z),imag(z),real(w),imag(w)); 
     dZw = min(abs(Z-w.'),[],2);  % distance to nearest corner
     ii = find(dZw>5e-3); outside = ~inpolygonc(zz,Z(ii));
     uu = abs(uv(zz)); uu(outside) = NaN; umax = max(max(uu));
     pcolor(x,y,uu), hold on, colormap(gca,parula)
     shading interp, colorbar, caxis([0 umax])
     plot(Z([1:end 1]),'k',LW,.8)
     pp = psi(zz); pp(outside) = NaN; pmin = min(min(pp)); pmax = max(max(pp));
     lev = pmin+(.1:.1:.9)*(pmax-pmin);
     contour(x,y,pp,lev,'k',LW,.6)
     psiratio = pmin/pmax; fac = max(psiratio,1/psiratio);
     if sign(fac) == -1     % Moffatt eddies in yellow
        lev1 = lev(1:2:end)*fac;
        contour(x,y,pp,lev1,'y',LW,.55)
        if abs(fac) > 1e4   % second eddies (white)
           lev2 = lev1*fac; contour(x,y,pp,lev2,LW,.5,CO,.99*[1 1 1])
        end
        if abs(fac) > 1e3   % third eddies (yellow)
           lev3 = lev2*fac; contour(x,y,pp,lev3,'y',LW,.45)
        end
     end
     if nargin==5, plot(cell2mat(Pol),'.r',MS,8), end
     hold off, axis([xm+.6*dx*[-1 1] ym+.6*dy*[-1 1]]), axis equal off
\end{verbatim}
\par}
\caption{Contour plotter.}
\end{figure}

\section*{Acknowledgments}
We have benefitted from helpful advice from Peter Baddoo,
Martin Bazant, Abi Gopal, Yuji Nakatsukasa, Hilary and John Ockendon, and
Kirill Serkh.

\end{document}